\magnification \magstep1



     \hsize 6.3 truein

     \font\twelverm=cmr10 scaled \magstep1
     \font\ninerm=cmr9
     \font\fifteenrm =cmr10  scaled \magstep2
     \font\numsetfont=msbm10

     \def\centered#1{\hfil \break \centerline{#1}}
     \def\newline{\hfil \break}

     \def\balid{\hbox{,\kern-.08em ,}}
     \def\endofproof{\ \vrule height 8pt depth 0pt width 7pt \medskip}

     \def\Re{\mathop {{\rm Re}} \nolimits}

     \def\numset#1{\hbox{\numsetfont #1}}

 \def\abstract#1{{\midinsert \narrower \ninerm {\bf Abstract.} #1 \endinsert}}

 \def\squar#1{{\vcenter{\vbox{\hrule height.4pt \hbox{\vrule width.4pt
height#1pt \kern#1pt \vrule width .4pt} \hrule height .4pt}}}}

     \def\text#1{\ \hbox{\rm #1}}

     \baselineskip=13pt plus 1pt minus 1pt \parskip=3pt

     \def\title#1{\null\bigskip \centered{{\fifteenrm #1 }}\bigskip }

     \def\author#1{{\noindent \twelverm #1}}

     \newcount\allszam \allszam=0
     \newcount\thszam \thszam=0

     \def\section#1#2{\vskip 30pt \goodbreak \centerline{{\twelverm #1\ #2}}\nobreak
\bigskip \def\secno{#1} \allszam=1 }

     \def\claimno{{\secno\the\allszam.\ }}

     \def\Def#1{\medskip\goodbreak {\bf Definition \claimno}  #1
\medskip \advance\allszam by1}

     \def\Theorema#1#2{\medskip\goodbreak \ifnum\thszam=0 {\bf \claimno Theorem} {\
(#1).}   {\sl #2.} \medskip \advance\allszam by1 \else {\bf Theorem\
\the\thszam.}{\ (#1).}   {\sl #2} \medskip \advance\thszam by1 \fi}
     \def\Theorem#1{\medskip\goodbreak \ifnum\thszam=0 {\bf Theorem \claimno}
     {\sl #1} \medskip \advance\allszam by1 \else {\bf Theorem\ \the\thszam.}
{\sl #1} \medskip \advance\thszam by1 \fi}

     \def\Lemma#1{\medskip\goodbreak {\bf Lemma \claimno}   {\sl #1} \medskip
\advance\allszam by1}

     \def\St#1{\medskip\goodbreak {\bf Statement \claimno}   {\sl #1}
\medskip \advance\allszam by1}

     \def\Prop#1{\medskip\goodbreak {\bf Proposition \claimno}   {\sl #1}
\medskip \advance\allszam by1}

     \def\Cor#1{\medskip\goodbreak {\bf Corollary \claimno}   {\sl #1}
\medskip \advance\allszam by1}

     \def\Rem#1{\medskip\goodbreak {\bf Remark \claimno} #1 \medskip
\advance\allszam by1}
     
     \def\Conj#1{\medskip\goodbreak {\bf Conjecture \claimno} #1 \medskip
\advance\allszam by1}

     \def\Prob#1{\medskip\goodbreak {\bf Problem \claimno} #1 \medskip
\advance\allszam by1}
     \def\Proof{{\it Proof\relax}}

        \def\EQU#1{\eqno{(\secno#1)}}
     \def\equ#1{(\secno#1)}
     \def\authorhead#1{}
     \def\titlefej#1{}

     \def\printed{Printed \number\year.\number\month.\number\day.}
     \thszam=1

     \title{Prime values of reducible polynomials, II.}

     \author{Yong-Gao Chen,\footnote{${}^1$}{\tenrm Supported by National
Natural Science Foundation of China, Grant No.~0171046 and the `` 333 Project"
Foundation of Jiangsu Province of China. The work was done while first author
was visiting the Mathematical Institute of the Hungarian Academy of Sciences.}
     G\'abor Kun,
     G\'abor Pete,\footnote{${}^2$}{\tenrm Supported by Hungarian National
Foundation for Scientific Research, Grants No. F 026049 and T 30074}
     Imre Z. Ruzsa\footnote{${}^3$}{\tenrm Supported by Hungarian National
Foundation for Scientific Research, Grants No. T 025617 and T 29759}
     and \'Ad\'am Tim\'ar}

     \authorhead{Chen, Kun, Pete, Ruzsa, Tim\'ar}

     \medskip  \noindent
     Department of Mathematics, Nanjing Normal University
     \hfil\break Nanjing 210097, Jiangsu, P. R. China
     \hfil\break e-mail:ygchen@pine.njnu.edu.cn
     \medskip           \noindent
     Department of Computer Science, E\"otv\"os University
     \hfil\break Budapest, Kecskem\'eti utca 10-12, H--1053 Hungary
     \hfil\break e-mail: kungabor@cs.elte.hu
     \medskip           \noindent
     Bolyai Institute, University of Szeged
     \hfil\break Szeged, Aradi v\'ertan\'uk tere 1., H--6720 Hungary.
     \hfil\break e-mail: gpete@sol.cc.u-szeged.hu
     \medskip \noindent
     Mathematical Institute of the Hungarian Academy of Sciences,
     \hfil\break Budapest, Pf. 127, H--1364 Hungary.
     \hfil\break e-mail: ruzsa@math-inst.hu
     \medskip           \noindent
     Bolyai Institute, University of Szeged
     \hfil\break Szeged, Aradi v\'ertan\'uk tere 1., H--6720 Hungary.
     \hfil\break e-mail: tadam@petra.hos.u-szeged.hu

     \bigskip

     \section{1.}{Introduction}

     It is a generally accepted conjecture that an irreducible integer-valued
polynomial without a constant divisor assumes infinitely many prime values at
integers. On the other hand it is easy to see that for a reducible $f\in \numset{Q}[x]$
there are only finitely many integers $n$ for which $f(n)$ is prime. It is,
however, a nontrivial question to estimate the number of these integers. We
shall be primarily interested in finding estimates in terms of the degree of
$f$.

     In the sequel by ``polynomial'' we always mean a polynomial with rational
coefficients, and reducibility is meant in $\numset{Q}[x]$. We will write
     $$ P(f) = \# \{m\in \numset{Z}: f(m) \text{is prime} \} . $$

     In the first part [3] we investigated the class of {\it integer-valued
polynomials}, that is, such polynomials that $f(n)$ is integral whenever so is
$n$. We proved that
     $$ \eqalign { \exp \left ( c {n\over  \log n } \right ) & <
     \sup \{ P(f): \deg f =n, f \text{ is integer-valued and reducible
in } \numset{Q}[x] \} \cr
     & < \exp \left ( C {n\over  \log n } \right ) } $$
     with positive absolute constants $c, C$.

     In this part we investigate the behaviour of $P(f)$ under further
restrictions. We shall assume that

     (a) the factors of $f$ are also integer-valued,

     or that

     (b)    $f$ has integral coefficients, in which case by Gauss' lemma we
may also assume that the factors have integral coefficients.

     These assumptions considerably reduce the possible number of prime values.
Indeed, if $f=gh$ with integer-valued $g$ and $h$,
then $f(x)$ can be a prime only if either $g(x) = \pm  1$ or $ h(x) = \pm 1$, which
immeditately gives $2n$ as an upper bound. (Ore [5] attributes this observation
to St\"ackel [8].) Our aim is to improve this bound.

     The more natural case (b) has been investigated by Ore [5]. His result
sounds essentially as follows. [He formulates it indirecly (a polynomial which
assumes more than ... prime values must be irreducible), and does not give the
construction for general $n$.]

     \Theorem{
     Let $f\in \numset{Z}[x]$ be a reducible polynomial of degree $n$. If $n\ne 4, 5$, then
$P(f)\le n+2$. On the other hand, for every $n$ there is a reducible $f\in \numset{Z}[x]$ such
that $P(f)\ge n+1$. If $n=4$ or 5, then the maximal possible value of $P(f)$ is
8.}

     We think that the upper bound gives the truth.

     \Conj{ For every $n$ there is a reducible $f\in \numset{Z}[x]$ such that $P(f)=n+2$. }

     In Section 2 we outline Ore's argument, describe the construction and to
support the conjecture we show that it follows from certain generally accepted
but hopeless conjectures about primes.

     In case (a) we can also  reduce the trivial upper bound, though we are far
from  a complete answer.

     \Theorem{
     There is a constant $c  <2$ such that we have
     $$ P(f) < c  n + o(n)         $$
     for every polynomial of degree $n$ which can be written as a product of
two integer-valued nonconstant polynomials. A possible value of this constant
is $ c = 1.8723406362...$, determined by the formula $c=1+1/t$, where $t$ is
the only real root of the equation
     $$ t  \bigl(  2 \log t + 1/2  \bigr)  = 2 \log 2 - 1/2 . $$ }

     Besides $P(f)$ we will also consider
     $$ P^+(f) = \# \{m\in \numset{Z}: f(m)>0 \text{is prime} \} . $$
     This is a less natural concept; however, the restriction to positive
primes will enable us to give an exact answer. In the case considered in Part I
the discrepancy between the lower and upper bound was so large that this
distinction did not matter.

     \Theorem{
     For every $n\ge 2$ and every polynomial of degree $n$ which can be written as
a product of two integer-valued nonconstant polynomials we have $P^+(f)\le n$. On
the other hand, there is a reducible $f\in \numset{Z}[x]$ for which $P^+(f)=n$,
consequently the maximum of $P^+$ in both cases (a) and (b) is exactly $n$.}

     \section{2.}{Polynomials with integer coefficients}

     The upper bounds stated in Theorem 1 are due to Ore [5]. We outline his
argument, since the source is not easily available and also it gives some
background to the construction and the related conjecture.

     For a polynomial $f$, write
     $$ E^+(f) = \# \{m\in \numset{Z}: f(m) = +1\}, \quad E^-(f) = \# \{m\in \numset{Z}: f(m) = -1\},$$
     $$ E(f) = E^+(f) + E^-(f) = \# \{m\in \numset{Z}:  |f(m) | =  1\}. $$

     In this section by polynomial we will always mean a polynomial with
integer coefficients.

     The starting point is the following result of Dorwart and Ore [4].

     \Lemma{
     If a polynomial of degree $n$ satisfies $E(f)>n$, then $n\le 3$ and $f$ is of
the form $f(x) = \pm  h_i(\pm x+a)$, where the polynomials $h_i$, $i=1, ..., 5$ are
listed below: }
     \settabs\+ \indent \indent
        & $ h_1(x) = x(x-1)(x-3)+1$ , \qquad & $n=$ &3, \qquad $E(f)=$& 4 \cr
     \+ & $ h_1(x) = x(x-1)(x-3)+1$ , \qquad & $n=$ &3, \qquad $E(f)=$& 4 \cr
     \+ & $ h_2(x) =  (x-1)(x-2)-1$ , \qquad &      &2,             & 4 \cr
     \+ & $ h_3(x) =  2x   (x-2)+1$ , \qquad &      &2,             & 3 \cr
     \+ & $ h_4(x) =  2x-1        $ , \qquad &      &1,             & 2 \cr
     \+ & $ h_5(x) =   x-1        $ , \qquad &      &1,             & 2.\cr

     This immediately yields that if a reducible polynomial $f=gh$ satisfies
$P(f)>n$, then at least one of $g,h$ is from the above list. Furthermore, we
see that $P(f)\le n+4$, and if $P(f)=n+3$ or $n+4$, then both factors come from
the list. For $n+4$ prime values, the  only possibility is
     $$ f(x) = \pm  h_2(\pm x+a) h_2(\pm x+b), $$
     and indeed we get 8 prime values for
     $$ f(x) =  \bigl(  1+x(x-3) \bigr)   \bigl(  1+ (x-4)(x-7) \bigr)  . $$
     For $n+3$ prime values one of the factors must be an $h_2$, and the other
factor has to assume prime values at 4 consecutive integers. However, $h_3$,
$h_4$ and $h_5$ are easily shown not to have this property by simple
divisibility arguments. Another factor of type $h_1$ is possible, and an
example is
     $$ f(x) =  \bigl(  1+x(x-3) \bigr)    \bigl(  1+ (x-4)(x-5)(x-7) \bigr) . $$

     This ends Ore's argument. We will now discuss the case of $n+1$ or $n+2$
prime values. The above facts show that for $n>6$ the only possibility to have
$n+2$ prime values happens when one factor is of type $h_2$.
     Thus we fix $h(x)=h_2(x) $. We show how to construct, for a given $n\ge 3$, a
polynomial $g$ of degree $n-2$ so that $f=gh$ assumes $n+2$ prime values. This
construction depends on two unproved conjectures. The first is that $h$ assumes
infinitely many primes; it is a generally accepted conjecture that this holds
for every irreducible polynomial without constant divisor. The next is that for
every finite collection $a_1, ..., a_k$ of nonzero integers we can find an
integer $t$ such that all the numbers $1+ta_i$ are prime, a version of the
prime $k$-tuple conjecture (we will use it for $k=4$).

     Assuming the first conjecture, let $b_1, ..., b_{n-2}$ be distinct
integers such that each $h(b_i)$ is a prime. We put
     $$ g(x) = 1 + t (x-b_1)...(x-b_{n-2}) $$
     with suitable $t$. For every choice of $t$ we have $f(b_i)=h(b_i)=$ prime.
Now the second conjecture yields the existence of a $t$ such that $g(i)$ is
prime for $i=0,1,2,3$, and then so is $f(i)$ since $h(i)=\pm 1$ for these numbers.

     Ore's arguments show that this is essentially the only choice of $h$,
hence the first conjecture is necessary.

     The second conjecture can possibly be weakened for our purposes. Indeed,
we do not need prime-yielding values of $t$ for every $b_1, ..., b_{n-2}$; what
we need is that from the set
     $B= \{b: h(b) \text{ is prime} \}$
     we can select some $n-2$ such that the prime-quadruple conjecture works
for the  four numbers, determined by these $n-2$ elements of $B$ in the above
described way. Hence an average version of the prime tuple conjecture,
similar to that proved by Balog [2], may suffice.

     Finally we prove unconditionally that $n+1$ prime values can be attained
for every $n\ge 6$. One of the factors must come from the list, and like in the
above conditional argument we need that it assume infinitely many primes. The
only polynomials for which this is established are the linear ones, thus we
have to use $h_4$ or $h_5$. We will use $h(x)=x=h_5(x+1)$.

     Let $p_1, ..., p_{n-1}$ be distinct (not necessarily positive) primes. We
put
     $$ g(x) = 1 + t (x-p_1) ... (x-p_{n-1}) $$
     with a suitable integer $t$. Then $f=gh$ satisfies $f(p_i)=p_i$ for $i=1,
.., n-1$ and
     $$ f(1) = g(1) = 1 + t(1-p_1) ... (1-p_{n-1}), $$
     $$ f(-1) = -g(-1) = -  \bigl( 1 + t(-1-p_1) ... (-1-p_{n-1})  \bigr) . $$

     In general it seems difficult to make two such expressions simultaneously
prime. We get around this difficulty by selecting distinct primes $p_1, ...,
p_{n-1}$ so that
     $$ (1-p_1) ... (1-p_{n-1}) = (-1-p_1) ... (-1-p_{n-1}). \EQU1 $$
     This will guarantee $g(-1) = g(1)$ independently of the choice of $t$, and
if we select $t$ to make these numbers prime, then $f(1)=g(1)$ and
$f(-1)=-g(1)$ will be prime besides $f(p_i)=p_i$. This can be done by
Dirichlet's theorem.

     To achieve \equ1, if $n$ is odd, we simply use primes in pairs with their
negatives, that is $p_2=-p_1$, $p_4=-p_3$ and so on. Every such pair
contributes the same to both products.

     If $n\ge 4$ is even, we use primes in pairs except the last three which will
be $2, -3$ and $-5$. These contribute $ (-1)\cdot 4\cdot 6 = (-3)\cdot  2 \cdot  4 = -24$ to both
sides.

     Finally we mention two examples that establish the maximal value for
degree 2 and 3:

     $n=2$: $f(x) =  \bigl( x \bigr)   \bigl( x-4 \bigr) $: $P(f)=4,$

     $n=3$: $f(x) =  \bigl(  1+x(x-3) \bigr)    \bigl( x-5 \bigr) $: $P(f)=5.$

     \section{3.}{Integer-valued polynomials}

     In this section we prove Theorem 2.

     \Lemma{
     Let $a_1, ..., a_k, b_1, ..., b_k$ be $2k$ distinct integers. Write
     $$ U = \prod _{i<j}  |a_i-a_j |, \quad V = \prod _{i<j}  |b_i-b_j |, \quad
     D = \prod _{i,j}  |a_i-b_j |. $$
     We have $D \ge  UV (4/9)^k$. }

     This is Lemma 2.1 of [6].

     \Lemma{
     With the above notations we have
     $$ D \ge  (2/3)^k  \bigl(  1! 2! ... (2k-1)!  \bigr) ^{1/2} . \EQU1 $$ }

     \Proof.
     Let $c_1<c_2 < ... < c_{2k} $ be the sequence of the integers $a_1, ...,
a_k, b_1, ..., b_k$ arranged in increasing order. Clearly
     $$ W = \prod _{i<j} (c_j-c_i) \ge  \prod _{1\le i<j\le 2k} (j-i) = 1! 2! ... (2k-1)! . $$
     On the other hand we have $W = UVD \le  D^2 (9/4)^k$ by the previous lemma.
     \equ1 follows by comparing these inequalities. \endofproof

     \Lemma{
     Let $a_1, ..., a_k, b_1, ..., b_s$ be $ k+s$ distinct integers, $k\le s$.
     We have
     $$ D = \prod _{i=1}^k \prod _{j=1}^s  |a_i-b_j | \ge  (2/3)^s  \bigl(  1! 2! ... (2k-1)!
 \bigr) ^{s/(2k)} . \EQU2 $$ }

     \Proof. By the previous lemma we have
     $$ D = \prod _{i=1}^k \prod _{j=1}^k  |a_i-b_{m_j} | \ge  (2/3)^k  \bigl(  1! 2! ... (2k-1)!
 \bigr) ^{1/2} $$
     for an arbitrary sequence $m_1, ..., m_k$ of distinct integers satisfying
$1\le m_j\le s$. Multiplying these inequalities for all possible choices of the $m_j$
and taking an appropriate root we
obtain \equ2. \endofproof

     \Proof\ of the Theorem.
     Let $f$ be an integer-valued polynomial of degree $n$. We shall find an
upper estimate for $E(f)$ in the form $cn + o(n)$ with the $c$ given in the
theorem.

     Write $r= E^+(f)$, $s= E^-(f) $
     and take integers $a_1, ..., a_r, b_1,..., b_s$ so that $f(a_i)=1$,
$f(b_j)=-1$. The polynomial $F=n!f$ has integer coefficients. Write
     $$ F(x) + n! = A \prod _{j=1}^n (x-\beta _j) . \EQU3 $$
     $A$ is an integer, hence $ |A |\ge 1$.
     Since the integers $b_j$ are roots of the polynomial $F(x)+n!$, they are
listed among the $\beta _j$, say $\beta _1=b_1,$, ..., $\beta _s=b_s$.

     We substitute $x=a_i$ into \equ3 to obtain
     $$ 2n! = F(a_i) + n! = A \prod _{j=1}^n (a_i-\beta _j) . \EQU4 $$

     For each $s+1\le j\le n$ there may be at most one $i$ for which
     $$ -1/2 \le  a_i - \Re \beta _j < 1/2 . \EQU5 $$
     This makes altogether at most $n-s$ values of $i$, so there are at least
$r-(n-s)=(r+s)-n$ values for which $a_i$ does not satisfy any of the
inequalities \equ5. We may assume that these are $a_1, ..., a_k$, where
$k=r+s-n$. We may also assume that $k>n/2$, since otherwise $E(f) = k+n \le  (3/2)
n $ and we are ready.

     Now we multiply the equations \equ4 for $i=1, ..., k$. This yields
     $$ (2n!)^k = A^k \prod _{i=1}^k \prod _{j=1}^n (a_i-\beta _j)
                = A^k \prod _{i=1}^k \prod _{j=1}^s (a_i-b_j)
                      \prod _{i=1}^k \prod _{j=s+1}^n (a_i-\beta _j). \EQU6 $$
     Observe that $k=r+s-n\le s$.

     Now we give a lower estimate for the right side of \equ6. Take first a $j$
satisfying $s+1\le j\le n$. We have
     $$  |a_i - \beta _j | \ge   | a_i - \Re \beta _j  | . $$

     If we arrange all the numbers $  |m - \Re \beta _j |$, $m\in \numset{Z}$ in increasing order,
we get the sequence
      $ \gamma , 1-\gamma , 1+\gamma , 2-\gamma , 2+\gamma , 3-\gamma , ..., $,
      where $\gamma $ is the distance of $ \Re \beta _j$ from the nearest integer. The
factors of our product are $  k$ numbers from this sequence, and the first term
($\gamma $) is excluded. Since $1-\gamma \ge 1/2$, $1+\gamma \ge 1$, $2-\gamma \ge 3/2$ and so on, the product
of $  k$ terms is at least
     $$ {1\over 2} {2\over 2} {3\over 2} ... {  k\over 2} = {   k ! \over  2^k  } . $$
     In particular,
     $$ \prod _{i=1}^k  |a_i-\beta _j | \ge  2^{-k} k!, $$
     $$ \prod _{i=1}^k \prod _{j=s+1}^n  |a_i-\beta _j | \ge  2^{-k(n-s)} k!^{n-s} . $$

     To estimate the first double product in \equ6 we use Lemma 3.3, and we use
$ |A |^k\ge 1$. These inequalities together give
     $$ (2n!)^k \ge  (2/3)^s  \bigl(  1! 2! ... (2k-1)!  \bigr) ^{s\over  2k } 2^{-k(n-s)} k!^{n-s}
. $$
     By the symmetric role of $r$ and $s$ we also have
     $$ (2n!)^k \ge  (2/3)^r  \bigl(  1! 2! ... (2k-1)!  \bigr) ^{r\over  2k } 2^{-k(n-r)} k!^{n-r}
. $$
     We multiply these inequalities and we obtain (recall that $r+s=n+k$)
     $$ (2n!)^{2k} \ge  (2/3)^{n+k}  \bigl(  1! 2! ... (2k-1)!  \bigr) ^{{n+k}\over  2k }
2^{-k(n-k)} k!^{n-k} . $$

     To utilize this inequality we take logarithm of both sides and use the
familiar estimate
     $$ \log m! = m ( \log m -1) + O( \log m) $$
     and the following which can be deduced from it by an immediate calculation
     $$ \log  ( 1! 2! ... m! ) = m^2 \left ( {1\over 2} \log m - {3\over 4} \right ) + O(m \log m). $$
     We obtain
     $$ \eqalign {
     & 2kn ( \log n -1) \cr
     \ge  & k (n+k) ( \log k + \log 2 - 3/2) -k(n-k) \log 2 +k(n-k)
( \log k -1) + O(n \log n) . \cr } $$
     After dividing by $k^2$ and cancelling certain terms this inequality
becomes
     $$ {n\over k} \left ( 2 \log { n\over k} + {1\over 2} \right ) \ge  2 \log 2 - {1\over 2} + O\left ( \log n \over  n \right )
     . $$
     Thus we get $n/k\ge t + o(1)$, where $t$ is the solution of
     $ t  \bigl(  2 \log t + 1/2  \bigr)  = 2 \log 2 - 1/2 $. We find $t=1.1463411865...$
which leads to $r+s=k+n\le n(1+1/t) +o(n)$. The constant appearing here is
$ 1+1/t = 1.8723406362...$ . \endofproof

     {\bf Remark.}
     Let $S$ be an arbitrary finite set. Let $E_{S}(f)$ denote the number of
distinct integers $a$ such that $f(a) \in S$. The proofs of Theorems 1 and 2
depended on estimations for the value of $E_{S}(f)$ in the case that $S =
\{-1,1\}$. In the case of more general finite sets $S$ similar estimates can be
obtained.

     Following the approach of Dorwart and Ore, one can show that if $f$ has
integer coefficients, then $E_{S}(f) \le n$ except for a finite list of
polynomials and their translations; in particular, $E_{S}(f) \le n$ for $n$
sufficiently large. However, it seems to be a nontrivial question to find sharp
estimates (in terms of the set $S$) for the number and maximal degree of the
exceptional polynomials.

     For integer-valued polynomials $f$ we can show that $E_{S}(f) < C n +
o(n)$ for some absolute constant $C$ (independent of the size of the finite set
$S$). This can be done by modifying the proof of Theorem 2, and in this way we
obtained $C=3$.

     We indicate a different proof that yields a somewhat better constant. Let
$K$ denote the maximum of absolute values of elements of $S$. A theorem of
P\'olya [6] (see also in Aigner-Ziegler [1]) asserts that for a polynomial $f$ of
degree $n$ and leading coefficient 1 the measure of real numbers satisfying
$ |f(x) |\le 1$ is at most 4 (in fact, the measure of the real parts of such complex
numbers $x$ is at most 4). By a natural rescaling we get that if the leading
coefficient is $c$, then the measure of reals satisfying $ |f(x) |\le K$ is at most
$4(K/ |c |)^ {1/n}$. Since this set is the union of at most $n$ intervals, the
number of integers satisfying $ |f(k) |\le K$ is at most
     $ n + 4(K/ |c |)^ {1/n}$.
     Since for integer-valued polynomials we have $ |c |\ge 1/n!$, we obtain
     $$ E_S(f) \le  n + 4(Kn!  )^ {1/n} = \left ( 1 + {4\over e} \right ) n + O( \log n). $$

     On the other hand, we do not have any better lower bound than $E_S(f)\ge
n+2$ for even values of $n$, which follows by considering
     $$ f(x) = a {x \choose n} + b $$
     for suitable $a,b$. We cannot even achieve this for general odd $n$.

     \section{4.}{The case of positive primes}

     In this section we prove Theorem 3.

     To show that $n  $ prime values are possible we apply the following
construction. Let $p_1, ..., p_{n-1}$ be distinct positive primes. We put
$h(x)=x$ and
     $$ g(x) = 1 + t (x-p_1) ... (x-p_{n-1}) $$
     with a suitable integer $t$. Then $f=gh$ satisfies $f(p_i)=p_i$ for $i=1,
.., n-1$ and
     $$ f(1) = g(1) = 1 + t(1-p_1) ... (1-p_{n-1}). $$
     This will be a positive prime for a suitable choice of $t$ by Dirichlet's
theorem.

     Next we show that the number of prime values is at most $n$.

     Let $f=gh$, where $g,h$ are integer-valued polynomials of degree at  least
1. If $f(m)$ is prime, then either $g(m)=\pm 1$ or $h(m)=\pm 1$. Hence the upper
estimate follows from the following statement.

     \St{
     Let $g, h$ be polynomials of degree at least one with real coefficients
and write $f=gh$. Consider those real numbers that satisfy
     \item{(a)}{ $ g(x) =\pm 1$ or $ h(x)=\pm 1$, and}
     \item{(b)}{ $f(x)>1 $.}
     \item
     The total number of such reals is at most $n = \deg f$.}

     \Proof. Let these numbers be $x_1< ...< x_k$. We will show that
     $$ \text{ (number of roots of $g'$) + (number of roots of $h'$)} \ge  k-2.
\EQU1 $$
     This clearly implies the statement.

     We divide the points $x_i$ into four types. It is of type $g+$ if
$g(x_i)=1$; the types $g-$, $h+$ and $h-$ are defined analogously. By a {\it
block} we mean a maximal sequence of consecutive $x_i$'s of the same type; the
type of the block is this type. Let $l$ denote the number of blocks. The number
of pairs $x_i, x_{i+1}$ of equal type is then exactly $k-l$.

     We call a block {\it extremal}, if it contains $x_1$ or $x_k$, and {\it
central} otherwise. The number of extremal blocks is 1 or 2, the number of
central blocks is at least $l-2$.

     We will show that
     $$ \eqalign {
     \text{ number of roots of } g' \ge  & \text{ (number of pairs of type } g\pm )
\cr + & \text{ (number of central blocks of type }h\pm  ), }\EQU2 $$
     and similarly
     $$ \eqalign {
     \text{ number of roots of } h' \ge  & \text{ (number of pairs of type } h\pm )
\cr + & \text{ (number of central blocks of type }g\pm  ). }\EQU3 $$
     On adding these inequalities we get the left side of \equ1, and on the
right side we have at least $(k-l)+(l-2)=k-2$ as claimed.

     To prove \equ2 we are going to map the pairs of type $g\pm $ and blocks of
type $h\pm $ onto roots of $g'$.

     Given a pair of type $g\pm $ we have $g(x_i)=g(x_{i+1})$. Hence $g'$ has at
least one root in the interval $(x_i, x_{i+1})$. We map this pair onto this
root (or onto any of these roots, if there are more than one).

     Consider now a central block of type $h\pm $, say $(x_i, ..., x_j)$, where
$1<i\le j<k$ by definition. For sake of definiteness assume it is of type $h+$. We
have $h(x_i) = ... = h(x_j)=1$, so $g(x_i)> 1, ..., g(x_j)> 1$. On the other
hand, $x_{i-1}$ and $x_{j+1}$ are of a different type, consequently
$g(x_{i-1})\le  1$, $g(x_{j+1})\le  1$. We map this block onto any local place of
maximum of $g$ within the interval $[x_{i-1}, x_{j+1}]$. The previous
inequalities show that this cannot be any of the endpoints, thus it must be a
root of $g'$.

     We are going to show that we use any given root at most once. This has
three subcases.

     The roots corresponding to pairs are obviously distinct.

     A root corresponding to a block is within an interval $[x_t, x_{t+1}]$,
where at least one of $x_t$ and $x_{t+1}$ is a member of that block, in
particular, it is of type $h\pm $. This shows that it cannot coincide with a root
corresponding to a pair.

     Consider finally two blocks, say $(x_i, ..., x_j)$ and $(x_u, ..., x_v)$,
such that $j<u$. The corresponding roots are situated in the intervals
$(x_{i-1}, x_{j+1})$ and $(x_{u-1}, x_{v+1})$, respectively. These are disjoint
unless $u=j+1$. If $u=j+1$, then the two blocks are adjacent, hence they must
be of different types, one of type $h+$ and the other of type $h-$. Hence $g$
has a local maximum at  one and a local minimum at the other, so they are
distinct.

     The proof of \equ3 proceeds in the same way, with the roles of $g$ and $h$
interchanged. \endofproof

     {\bf Remark.}
     Learning this result, N\'andor Sim\'anyi pointed out that Statement 4.1
holds for an arbitrary ordered field. The above proof can probably be extended
for this more general case; however, we can also argue as follows. For a fixed
pair of degrees $(\deg g,\deg h)$ this is a first order formula in the theory
of really closed ordered fields. This theory is complete, and we already know
that the statement is true for $\numset{R}$, therefore it is true for an
arbitrary really closed field. Finally, every ordered field has a really closed
extension, and the validity of the statement descends to subfields.

     He also asked whether a generalization of Statement 4.1 could be
valid for the case of complex polynomials $g,h$ and complex values of
$x_1,\dots,x_k$. We found that the answer is ``no", as shown
by the following example: $$g(z)={1\over3}z^3-z+1,\qquad
h(z)={2\over9}(z-2)^2+1.$$ Here $\deg(gh)=5$ but we have 6 ``bad" $x_i$'s:
$g(0)=g(\pm\sqrt3) =h(2)=1$, $h(2\pm 3i)=-1$, and $gh(x_i)\in (1,\infty)$
for all of them.

     This example now raises the question whether the maximal number of
possible complex $x_i$'s is equal to the trivial upper bound $2(\deg g+\deg
h)$. So far we have not been able to find more than 6 ``bad'' values for the
pair of degrees $(2,3)$; a reason for this can be B\'ezout's theorem on the
number of intersections of the real algebraic curves $g^{-1}(\numset{R} )$ and $h^{-1}(\numset{R}
)$.

     {\bf Acknowledgement.} We are grateful to a referee for several
corrections and for suggesting the problem of estimating $E_S(f)$ for general
sets (the Remark at the end of Section 3), and to John Rickard for the
information on P\'olya's inequality.

     \section{}{References}


      \newcount\refno \refno=0

     \def\paper#1#2#3#4#5#6{\advance\refno by1 \item{[\the\refno]}{#1, #3,
{\it #4} {\bf #5} (#2), #6.}}
     \def\toappear#1#2#3{\advance\refno by1 \item{[\the\refno]}{#1, #2, {\it #3,}
to appear.}}
     \def\cisu#1#2#3{\advance\refno by1 \item{[\the\refno]}{#1, #2, submitted
to the {\it #3.} }}
     \def\preprint#1#2{\advance\refno by1 \item{[\the\refno]}{#1, #2, preprint.}}
     \def\extra#1#2{\advance\refno by1 \item{[\the\refno]}{#1, #2.}}
      \def\book#1#2#3#4{\advance\refno by1 \item{[\number\refno]}{#1, {\it
 #3}, #4 #2.}}
      \def\proceedings#1#2#3#4#5#6{\advance\refno by1 \item{[\the\refno]}{#1, #3,
 {in: \it #4}, {#5} (#2), #6.}}
      \def\resx#1#2#3#4#5{\advance\refno by1 \item{[\the\refno]}{#1, #3,
 {in: \it #5},  #4 }}
      \def\resy#1#2#3#4#5#6{\advance\refno by1 \item{[\the\refno]}{#1, #3,
 {in: \it #4, #6},  #5 }}

     \frenchspacing

     \book {Aigner, M.; Ziegler, G. M} {1999} {Proofs from the book} {Springer}

     \proceedings {Balog, A.} {1990} {The prime $k$-tuplets conjecture on average}
{Analytic number theory (Allerton Park 1989)} {Birkh\"auser, Boston, MA}
{47--75}

     \paper {Chen, Y. G., Ruzsa, I. Z.} {2000} {Prime values of reducible
polynomials I.} {Acta Arithmetica} {95} {185--193}

     \paper {Dorwart, H.  L., Ore, O.} {1933} {Criteria for the irreducibility
of polynomials} {Annals of Math.} {34} {81--94}

     \paper {Ore, O.} {1934} {Einige Bemerkungen \"uber Irreduzibilit\"at}
{Jahresber. d. Deutsch. Math. Verein.} {44} {147--151}

     \paper {P\'olya, Gy.} {1928} {Beitrag zur Verallgemeinerung des
Verzerrungssatzes auf mehrfach zu\-sammenh\"agenden Gebieten} {Sitzungsber.
Preuss. Akad. Wiss. Berlin} {} {228-232; also in Collected papers, Vol. 1., MIT
Press 1974, 347--351}

     \paper {Ruzsa, I. Z.} {1992} {Large prime factors of sums} {Studia Sci.
Math. Hungar.} {27} {463--470}

     \paper {St\"ackel, P.} {1918} {Arithmetische Eigenschaften ganzer
Funktionen} {Journ. f. Math.} {148} {101--112}

     \printed \bye